.
\font\piccolo=cmr8.
\font\script=eusm10.
\font\sets=msbm10.
\font\stampatello=cmcsc10.
.
\def\0{{\bf 0}}
\def\1{{\bf 1}}
\def\defineq{\buildrel{def}\over{=}}
\def\EqByDef{\buildrel{\bullet}\over{=}}

\def\QED{\hfill {\rm QED}}
\def\C{\hbox{\sets C}}
\def\N{\hbox{\sets N}}
\def\sQ{\hbox{\script Q}}
\def\R{\hbox{\sets R}}
\def\Primes{\hbox{\sets P}}

\def\Z{\hbox{\sets Z}}
\def\square{\hbox{\vrule\vbox{\hrule\phantom{s}\hrule}\vrule}}
\def\supporto{{\rm supp}}

\def\WA{(\hbox{\stampatello WA})}

\def\BH{(\hbox{\stampatello BH})}
\def\GRE{(\hbox{\rm G.R.E.})}

\def\REEF{(\hbox{\stampatello R.E.E.F.})}
\def\RamaExp{(\hbox{\stampatello Ramanujan Expansion})}
\def\LuchtExp{(\hbox{\stampatello Lucht Expansion})}
\def\etaDecay{(\eta-\hbox{\stampatello Decay})}
\def\Cdp{(\hbox{\stampatello Condition}\;dp)}

\def\IrrP1F{{\hbox{\rm Irr}^{(P)}_1\,F}}

\def\Car{{\rm Car}}
\def\Win{{\rm Win}}
\def\CarT{\Car\; }
\def\WinT{\Win\; }

\par
\centerline{\bf Good Ramanujan Expansions:}
\par
\centerline{\bf A suitably enhanced decay of coefficients}
\par
\centerline{\bf has important consequences}
\bigskip
\centerline{Giovanni Coppola}\footnote{ }{MSC $2020$: $11{\rm N}37$ - Keywords: arithmetic function, Ramanujan expansion, coefficients decay} 

\bigskip

\par
\noindent
{\bf Abstract}. {\piccolo In this self-contained short note, we introduce the new definition of Good Ramanujan Expansion, say G.R.E., for a fixed arithmetic function $F$, building upon a good decay of its coefficients $G$; this, gains $\log-$powers w.r.t. the trivial bound for $G$ and precisely $\log^{1+\eta}$, where the present parameter $\eta>0$ is real. This property alone has important consequences for all the $F$ having a G.R.E. : mainly, 1) the Eratosthenes Transform $F'$ of our $F$ is infinitesimal (see in Theorem 1); 2) when $\eta>1$ (an enhanced decay) we have uniqueness of $G$ (actually, these are the classic Wintner-Carmichael coefficients, see Th.2); 3) we get a bound for $F$ (in Th.3); 4) an important new class of arithmetic functions $F$ can't have a G.R.E. (see Th.4). These are a generalization of Correlations; which in this way, if are, say, a kind of \lq \lq far from constants\rq \rq, may not have a G.R.E., whence, a fortiori, can't have the (R.E.E.F.). This is the Ramanujan Exact Explicit Formula, that we introduced with Prof. Ram Murty. 
On the Hardy-Ramanujan Journal, I proved that: any\thinspace \lq \lq {\stampatello fair}\rq \rq\thinspace Correlation may be \lq \lq well-approximated\rq \rq \thinspace by a (BH) Correlation. 
For a (BH) Correlation, we prove here: it has a G.R.E. $\Leftrightarrow$ 
it has the (R.E.E.F.). On the other hand, next Counterexample 1 in section 3 is a (BH) Correlation without the (R.E.E.F.); also, it provides here a (BH) Correlation behaving like in Theorem 4.}  

\bigskip

\par
\noindent{\bf 1. Introduction. Notations and New definitions} 
\bigskip
\par
\noindent
In 1918 Ramanujan [R] published a paper that is both a milestone in number theory and the first page of mathematical literature in the subject of {\stampatello Ramanujan Expansions for Arithmetic Functions}. We say that, once {\stampatello fixed} \enspace $F:\N \rightarrow \C$, we have {\stampatello a Ramanujan Expansion for } $F$, {\stampatello with coefficient} \enspace $G:\N \rightarrow \C$, IFF (if and only if) {\stampatello by definition}:
$$
\forall a\in \N,
\quad
F(a)=\lim_x \sum_{q\le x}G(q)c_q(a)\defineq \sum_{q=1}^{\infty}G(q)c_q(a), 
$$ 
\par
\noindent
where it's {\stampatello implicit}: ${\displaystyle \lim_x}$ {\stampatello is over} $x\to \infty$ {\stampatello and it exists in complex numbers}. 
\par
Here {\stampatello we abbreviate with} \enspace ${\displaystyle c_q(a)\defineq \sum_{{j\le q}\atop {(j,q)=1}}\cos {{2\pi ja}\over q} }$ \enspace {\stampatello the Ramanujan Sum of modulus $q\in \N$ and argument $a\in \N$}. (In this paper, {\stampatello finite \& infinite sums are over naturals : within} $\N$.) 
\smallskip
\par
We'll not give all the {\stampatello details about notations}, since we'll follow the standard ones; like for all {\stampatello implicit conventions \& standard facts, we follow $\S1$ of [C1]}. 
\smallskip
\par
The {\stampatello idea of Ramanujan} was (and still is!) {\stampatello simple and far-reaching : any $F:\N \rightarrow \C$ has a Ramanujan Expansion (see Remarks in $\S4$)}; He soon discovered that $\underline{\hbox{\stampatello it's not unique}}$ , {\stampatello whatever is our} $F$. This {\stampatello follows from expanding the null-function}, write $\0$, which sends \enspace $\N$ \enspace in the set $\{0\}$ : compare $\S1$ of [C1]. Ramanujan gave [R] explicit non-zero Coefficients for it, see $R_0(q)$, after Proposition 2. 
\par
In this short note, {\stampatello  we} {\stampatello prove how assuming a suitable decay in the}, say, {\stampatello $q-$th coefficients \enspace $G(q)$ \enspace in expansion} above {\stampatello entails important consequences}:
\item{1)} {\stampatello the Eratosthenes Transform} \enspace $F'\defineq F \ast \mu$, \enspace see [C1], of our $F$ {\stampatello is infinitesimal};  actually, we'll {\stampatello derive a bound} for $F'$ strictly {\stampatello depending on $G(q)$ decay}: see Theorem 1; whence Corollary 1 \& 2;
\item{2)} {\stampatello when the decay is \lq \lq enhanced\rq \rq}, say, w.r.t.(with respect to) our standard decay assumption, then {\stampatello we get} the $\underline{\hbox{\stampatello uniqueness of coefficients}}$ : details in Theorem 2; whence, Proposition 2; 
\item{3)} a {\stampatello Good Ramanujan Expansion}, abbrev. $\GRE$ , see our new definition, upcoming next, {\stampatello entails also a bound for our $F$ that we expand} : precise bound's in Theorem 3;    
\item{4)} a particular sub-class of periodic arithmetic functions, that we introduce in next definition, has {\stampatello members having Eratosthenes Transform}, say $F'$, with a particular \& characteristic behavior, ensuring $F'$ is {\stampatello not infinitesimal, excluding for them any $\GRE$}, by Theorem 1; and this is our last main result in this paper : Theorem 4. 

\bigskip

\par				
\centerline{\bf Notations: recalling standard \& non-standard}

\medskip

\par
\noindent
For standard notations, like numbers sets (by the way, $\Primes$ is the primes' set), we follow [C1] : there we find M\"obius Function $\mu$ \& Euler Totient Function $\varphi$ (with basic properties). Here, for the non-standard ones, namely those that may differ from one Author to the other, we recall $\log$ is in natural base $e$; we'll call $q-$periodic a function with minimal {\it period} {\tt T} dividing $q$, but may be {\tt T}$<q$ (see [C1]); notation $\Z_q^*$ is standard for the {\it reduced residue classes} modulo $q\in \N$, but we abuse it writing, for $q-$periodic functions' argument $a$, {\stampatello often $a\in \Z_q^*$ to mean}: $a\in \N$, $a\le q$, $(a,q)=1$. We write $p\equiv r(q)$ for $p\equiv r(\!\bmod \,q)$. For fixed $F:\N \rightarrow \C$  
$$
\Win_q\, F\defineq \sum_{d\equiv 0(\!\!\bmod q)}{{F'(d)}\over d} , 
$$  
\par
\noindent
whenever the series converges pointwise in $\C$, is the {\stampatello $q-$th Wintner coefficient} for $F$; varying $q\in \N$, we get {\stampatello Wintner's Transform} $\WinT F$: $q\in \N \mapsto \Win_q\, F\in \C$, a new arithmetic function that exists provided all $\Win_q\, F$ exist. Analogously, 
$$
\Car_q\, F\defineq {1\over {\varphi(q)}}\lim_x {1\over x}\sum_{a\le x}F(a)c_q(a) , 
$$  
\par
\noindent
whenever the limit exists in $\C$, is the {\stampatello $q-$th Carmichael coefficient} for $F$; varying $q\in \N$, we get {\stampatello Carmichael's Transform} $\CarT F$: $q\in \N \mapsto \Car_q\, F\in \C$, a new arithmetic function that exists provided all $\Car_q\, F$ exist.\hfill Just to abbreviate a limit over a subsequence of primes, that we need in the following, 
$$
\lim_{p\equiv r(\!\! \bmod q)}H(p)
$$
\par
\noindent
indicates the limit of $H(p)$, for given $H:\N \rightarrow \C$, over primes\enspace $p\equiv r(q)$, for fixed $q\in \N$ and (trivially) $r\in \Z_q^*$. 
\par
\noindent
Recall: once fixed \thinspace $r\in \N$, $F:\N \rightarrow \C$ \thinspace is \thinspace  {\stampatello $r-$even, by definition,} IFF we have\enspace $F(a)=F((a,r))$, $\forall a\in \N$. The {\stampatello cases $r=1$ and $r=2$ are trivial : } both entail that {\stampatello \lq \lq $r-$periodic\rq \rq \thinspace is equivalent to \lq \lq $r-$even\rq \rq.}
\par
\noindent
Like in [C1],\enspace $\EqByDef$\enspace means \lq \lq equals, from the definition\rq \rq;\thinspace $\QED$ is the \lq \lq end of a part\rq \rq \thinspace of a Proof;\enspace whose end's \hfill $\square$

\medskip

\centerline{\bf New Definitions}

\medskip

\par
\noindent
{\stampatello A Good Ramanujan Expansion}, abbrev. $\GRE$, for a fixed $F:\N \rightarrow \C$, means by definition that we have fulfilled the three requirements: 
$$
\forall a\in \N, 
\quad
F(a)=\lim_x \sum_{q\le x}G(q)c_q(a)
 \defineq \sum_{q=1}^{\infty}G(q)c_q(a)
\leqno{\RamaExp}
$$ 
\par
\noindent
where the {\stampatello Ramanujan Coefficient} $G:\N \rightarrow \C$ in this expansion is the same as in the
$$
\forall d\in \N, 
\quad
F'(d)=d\;\lim_x \sum_{K\le x}\mu(K)G(dK)
 \defineq d\,\sum_{K=1}^{\infty}\mu(K)G(dK)
\leqno{\LuchtExp}
$$ 
\par
\noindent
and it is infinitesimal, with the existence of a real {\stampatello parameter} $\eta>0$, in the following bound, valid $\forall q\ge 2$ : 
$$
G(q)\ll_{\eta} {1\over {q(\log q)^{1+\eta}}}
\leqno{\etaDecay}
$$ 
\par
\noindent
A G.R.E., in general, has only the requirement $\eta>0$ and may have {\it no Coefficients uniqueness}; however, compare Theorem 2. The Lucht Expansion is a vital hypothesis, here (see $\S4$ Remarks). 

\medskip

\par
A $\sQ-$periodic function \thinspace $F:\N \rightarrow \C$, with \thinspace $\sQ>2$, by definition, {\stampatello diverts values in $\Z_{\sQ}^*$} , \enspace when we have:\enspace $\exists a\in \Z_{\sQ}^*\backslash \{1\}$ with $F(a)\neq F(1)$. The other possibility, i.e. \thinspace $F(\Z_{\sQ}^*)=\{F(1)\}$, defines the $\sQ-$periodic functions $F$ that are {\stampatello $\sQ-$Monochromatic}. We may abbreviate, like in abstract, \lq \lq {\stampatello diverts values}\rq \rq \thinspace with \lq \lq {\stampatello far from constants}\rq \rq, since {\stampatello $\sQ-$Monochromatic} sounds like \lq \lq close to constant\rq \rq. 
\medskip
\par
\centerline{\bf Paper's Plan}
\smallskip
\item{$\star$} in next section we state and prove our results; 
\item{$\star$} in $\S3$ we apply some of our results to the Correlations of two arithmetic functions; 
\item{$\star$} in $\S4$ we offer some more Remarks \& properties, giving a glance to future applications for our results. 

\vfill
\eject

\par				
\noindent{\bf 2. Statements and Proofs of our Results} 
\bigskip
\par
A Good Ramanujan Expansion for $F$ entails an infinitesimal $F'$ ; more precisely, we have the following.
\smallskip
\par
\noindent{\bf Theorem 1.} {\it Let}\thinspace \thinspace $F:\N \rightarrow \C$ {\it have a}\thinspace \thinspace  {\stampatello Good Ramanujan Expansion}, {\it of}\thinspace \thinspace {\stampatello parameter } $\eta>0$. {\it Then}, $\forall d\ge 2$, 
$$
F'(d)\ll_{\eta} (\log d)^{-\eta}. 
$$
\par
\noindent
{\it In particular, the} {\stampatello Eratosthenes Transform $F'$ is infinitesimal}. 
\smallskip
\par
\noindent{\bf Proof}. Above $\LuchtExp$ with $G:\N \rightarrow \C$ of $\etaDecay$ gives, together with the {\stampatello easy inequality} : $\log(dK)\ge \max(\log d, \log K)$, $\forall d,K\in \N$, with \lq \lq {\stampatello large}\rq \rq \thinspace $d\in \N$ (actually, $d\ge 2$ suffices)
$$
F'(d)\ll_{\eta} \sum_{K\le d}{1\over {K(\log d)^{1+\eta}}}+\sum_{K>d}{1\over {K(\log K)^{1+\eta}}}
      \ll_{\eta} {1\over {\log^{\eta} d}}, 
$$
\par
\noindent
from the two easy bounds: 
$$
\sum_{K\le d}{1\over K}\ll \log d
\qquad
\hbox{\stampatello and}
\qquad
\sum_{K>d}{1\over {K(\log K)^{1+\eta}}}\le \int_{d}^{\infty}{1\over {t(\log t)^{1+\eta}}}dt
 \ll_{\eta} (\log d)^{-\eta}, 
$$
\par
\noindent
compare [T] for the former.\hfill $\square$ 
\medskip
\par
The {\stampatello divisor function}: \enspace $a\in \N \mapsto d(a)\defineq {\displaystyle \sum_{d|a}}\;1 \in \C$ \enspace {\stampatello can't have} a G.R.E., {\stampatello having Eratosthenes Transform} \enspace $\1(d)\defineq 1$, $\forall d\in \N$ \enspace {\stampatello not infinitesimal}. As a kind of {\stampatello curiosity, Ramanujan} [R] {\stampatello gave} 
$$
d(a)=\sum_{q=1}^{\infty}\left(-{{\log q}\over q}\right)c_q(a),
\quad \forall a\in \N
$$
\par
\noindent 
and, {\stampatello as you see}, these {\stampatello coefficients \lq \lq loose a $\log$\rq \rq \thinspace instead of gaining $\log^{1+\eta}$}. 

\medskip

\par
As a straight consequence of Theorem 1, we get
\smallskip
\par
\noindent{\bf Corollary 1.} {\it Let } $F:\N \rightarrow \C$ {\it have a } {\stampatello Good Ramanujan Expansion}, {\it of } {\stampatello parameter } $\eta>0$. {\it Then}
\smallskip
\centerline{
{\it as} \enspace $p\to \infty$ \enspace ({\it in prime numbers}),
\quad
${\displaystyle
F(2p)=F(2)+O_{\eta}\left((\log p)^{-\eta}\right)
 }$ 
}
\smallskip
\par
\noindent 
{\it and}, {\stampatello in particular}, \enspace ${\displaystyle \lim_p} \; F(2p)=F(2)$. {\it Furthermore, in the same way}, {\stampatello fixing} \enspace $a_0\in \N$,
\smallskip
\centerline{
{\it as} \enspace $p\to \infty$,
\quad
${\displaystyle
F(a_0 p)=F(a_0)+O_{\eta}\left((\log p)^{-\eta}\right)
 }$ 
}
\smallskip
\par
\noindent 
{\it whence}, {\stampatello in particular}, \enspace ${\displaystyle \lim_p} \; F(a_0 p)=F(a_0)$.  
\smallskip
\par
\noindent{\bf Remark 1.} The result may be {\stampatello generalized}: take {\stampatello a finite number of primes} \enspace $p_1<p_2<\cdots<p_r$ \enspace and then send \enspace $p_1\to \infty$, all staying in prime numbers, {\stampatello to multiply a fixed $a_0\in \N$}.\hfill $\diamond$

\medskip

\par
\noindent{\bf Proof}. Writing \enspace $F(a)=\sum_{d|a}F'(d)$, we get:\enspace $p>2\Rightarrow $\enspace $F(2p)-F(2)=F'(p)+F'(2p)$, whence Theorem 1 entails {\stampatello both} \enspace $F'(p)\ll_{\eta} (\log p)^{-\eta}$ \enspace {\stampatello and} \enspace $F'(2p)\ll_{\eta} (\log (2p))^{-\eta}\ll_{\eta} (\log p)^{-\eta}$.\hfill QED
\par
\noindent
Also,\enspace $p>a_0\Rightarrow $\enspace $F(a_0 p)-F(a_0)=\sum_{d|a_0}F'(dp)\ll_{\eta} d(a_0)\cdot (\log p)^{-\eta}$\enspace and\enspace {\stampatello $a_0\in \N$ fixed \thinspace $\Rightarrow$ \thinspace $d(a_0)$\thinspace fixed}.\hfill $\square$ 

\medskip

See that present Corollary 2, following, is \lq \lq much stronger\rq \rq, say, than Corollary 1. Our source of inspiration has been last line above in the Proof of Corollary 1. 
\vfill

\par
\noindent{\bf Remark 2.} The behavior of our $F$ in Corollary 1 resembles, also for the \lq \lq general\rq \rq \thinspace case $a_0\in \N$ and we assume it now {\stampatello to be even}, say $a_0=2k$, that of {\stampatello the expected heuristic for $F(2k)=$\lq \lq $2k-$twin primes count\rq \rq}, as given in {\stampatello Hardy-Littlewood Conjecture} (see [HL] Conjecture B). \enspace {\stampatello Focus on the \lq \lq Singular Series\rq \rq \thinspace : it has factors as a kind of finite products} \thinspace $\prod_{p|(2k)}\left(1+O(1/p)\right)$ \thinspace {\stampatello which, joining large factors $p$, say, \lq \lq don't change much\rq \rq}.\hfill $\diamond$

\eject

\par				
An even more important consequence of Theorem 1 is the following. 
\smallskip
\par
\noindent{\bf Corollary 2.} {\it Let } $F:\N \rightarrow \C$ {\it have a} {\stampatello (GRE)} {\it and be } $r-${\stampatello periodic}, {\it with} $r\in \N$ \& $r>2$. {\it Then } $F$ {\it is } $r-${\stampatello even}, {\it whence it has the} {\stampatello finite Ramanujan expansion}
$$
F(a)=\sum_{q|r}\widehat{F}(q)c_q(a),
\quad
\forall a\in \N, 
\leqno{(\ast)}
$$
\par
\noindent
{\stampatello with Unique Ramanujan Coefficients}: $\;$ $\widehat{F}(q)\defineq \Win_q\, F=\Car_q\, F$, $\forall q\in \N$. 
\smallskip
\par
\noindent{\bf Proof}. Theorem 1 above ensures that \enspace ${\displaystyle \lim_d}\, F'(d)=0$, which we'll use soon. 
\par
\noindent
Once {\stampatello fixed} $a\in \N$, $d_0:=(a,r)$ {\stampatello is fixed, too, like} $c:={a\over {d_0}}\in \Z_{{r\over {d_0}}}^*$ , and {\stampatello Dirichlet Theorem} [T] {\stampatello entails}:
$$
p\equiv c\left( \!\!\bmod {r\over {d_0}}\right)
\thinspace \Rightarrow \thinspace
F(a)=F(d_0 c)\buildrel{\hbox{\stampatello d.t.}}\over{=\!=\!=}\lim_p F(d_0 p)=\lim_p \sum_{d|d_0}F'(dp)+F(d_0)=F(d_0)=F((a,r)). 
$$
\par
\rightline{QED}
\par
\noindent
Then, {\stampatello from Ramanujan Orthogonality}: $\1_{d|a}={\displaystyle {1\over d}\sum_{q|d}c_q(a) }$, $\forall a,d\in \N$, see p.22[M], compare $(1.1)$ [C2], 
$$
\forall a\in \N, F(a)=F((a,r))
 =\sum_{{d|r}\atop {d|a}}F'(d)\buildrel{\hbox{\stampatello r.o.}}\over{=\!=\!=}\sum_{d|r}{{F'(d)}\over d}\sum_{q|d}c_q(a)
  =\sum_{q|r}\left( \sum_{{d|r}\atop {d\equiv 0(\!\!\bmod q)}}{{F'(d)}\over d} \right)c_q(a), 
$$
\par
\noindent
{\stampatello entailing $(\ast)$ above and $\widehat{F}(q)$ is well-defined, because uniformly $\forall q\in \N$}
$$
\sum_{{d|r}\atop {d\equiv 0(\!\!\bmod q)}}{{F'(d)}\over d}\EqByDef \Win_q\, F
\quad 
\hbox{\rm converges, since}
\quad 
\supporto(F')\subseteq \{d\in \N : d|r\}
$$
\par
\noindent
{\stampatello is finite}; {\stampatello for} the {\stampatello same reason}, we get {\stampatello Wintner Assumption}, see $\WA$ in next Theorem 2, where (in Proof's beginning) we also recall why $\WA$ entails\enspace $\Win_q\, F=\Car_q\, F$, $\forall q\in \N$.\hfill $\square$
\medskip
\par
\noindent 
Notice: the {\stampatello two limits} ${\displaystyle \lim_p }$ {\stampatello in the Proof above are}, since ${\displaystyle  p\equiv c\left( \!\!\bmod {r\over {d_0}}\right) }$,  {\stampatello actually} ${\displaystyle \lim_{{p\to \infty}\atop {p\equiv c\left( \!\!\bmod r/d_0\right)}} }$. 

\medskip

\par
{\stampatello With exactly same Proof as above (starting from 2nd-line) we have the more general:}
\smallskip
\par
\noindent{\bf Proposition 1.} {\it Let } $F:\N \rightarrow \C$ {\it have an} {\stampatello infinitesimal } $F'$ {\it and be } $r-${\stampatello periodic}, {\it with} $r\in \N$ \& $r>2$. {\it Then } $F$ {\it is } $r-${\stampatello even and $(\ast)$ above holds, with above Unique Ramanujan Coefficients.}
\medskip
We explicitly remark that the uniqueness of Ramanujan Coefficients is obtained, in Corollary 2, assuming a $\GRE$ with $\eta>0$, {\bf together with} the additional hypothesis: $F$ is $r-$periodic. Without this hypothesis, we need, say, $\eta>1$: see next result. 
\medskip
\par
{\stampatello Gaining} a little bit more in the {\stampatello decay of coefficients entails their uniqueness}. 
\smallskip
\par
\noindent{\bf Theorem 2.} {\it Let } $F:\N \rightarrow \C$ {\it have a } $\GRE$ {\it with} {\stampatello Ramanujan Coefficient} $G:\N \rightarrow \C$ {\stampatello and parameter } $\eta>1$. {\it Then}
\item{$a)$} the {\stampatello Transforms $\WinT F$ and $\CarT F$ exist} ; 
\smallskip
\item{$b)$} $\CarT F = \WinT F$ ; 
\smallskip
\item{$c)$} $G = \CarT F = \WinT F$. 

\medskip

\par
\noindent{\bf Remark 3.} In other words, point $c)$ says that any $F:\N \rightarrow \C$ with a G.R.E. featuring $\eta>1$ can only have as {\stampatello Ramanujan Coefficients} the so-called {\stampatello Carmichael-Wintner Coefficients}. \hfill $\diamond$ 

\medskip

\par				
\noindent{\bf Proof}. From \lq \lq {\stampatello Wintner's Criterion}\rq \rq, compare [C0] Proposition 2, point 1, above $a)$ and $b)$ follow from \lq \lq {\stampatello Wintner Assumption}\rq \rq (Criterion's hypothesis): 
$$
\sum_{d=1}^{\infty}{{|F'(d)|}\over d} < \infty.
\leqno{\WA}
$$
\par
\noindent
We apply, in order to prove $\WA$, the $\eta-${\stampatello decay} with $\eta>1$ to bound $F'(d)/d$, like in Theorem 1 Proof, getting 
$$
\sum_{d=1}^{\infty}{{|F'(d)|}\over d}\ll_{\eta} |F'(1)|+\sum_{d\ge 2}{1\over {d(\log d)^{\eta}}}
 < \infty, 
$$
\par
\noindent
from $\eta>1$.\hfill QED 
\par
\noindent
Using point $b)$, in order to get $c)$ we need to show $G=\WinT F$ : it comes from $(7)$ in Theorem 1 [C0], i.e. 
$$
\sum_{q=1}^{\infty}2^{\omega(q)}|G(q)|<\infty,
$$
\par
\noindent
recall \enspace $\omega(q)\defineq \left|\{p\in \Primes : p|q\}\right|$ \enspace is the number of prime factors of $q$; so this convergence may be proved as follows from the trivial $2^{\omega(q)}\le d(q)$, $\forall q\in \N$, whence
$$
\sum_{q=1}^{\infty}2^{\omega(q)}|G(q)|\ll_{\eta} |G(1)| + \sum_{q=2}^{\infty}{{d(q)}\over {q(\log q)^{1+\eta}}}
 \ll_{\eta} 1;
\leqno{(2.1)}
$$
\par
\noindent
in fact, {\stampatello Partial Summation} [T] {\stampatello entails} 
$$
\sum_{1<q\le x}{{d(q)}\over {q(\log q)^{1+\eta}}} = {1\over {x(\log x)^{1+\eta}}}\sum_{1<q\le x}d(q)
+\int_{2}^{x}\sum_{2\le q\le t}d(q)\,
{
{(t(\log t)^{1+\eta})'}
\over 
{t^2 (\log t)^{2+2\eta}}
}
dt
$$
\par
\noindent
and {\stampatello with classic bound} [T] :\enspace $\sum_{q\le x}d(q)\ll x\log x$, for all large $x\ge2$,\enspace we get
$$
\sum_{q=2}^{\infty}{{d(q)}\over {q(\log q)^{1+\eta}}}\EqByDef 
\lim_x \sum_{1<q\le x}{{d(q)}\over {q(\log q)^{1+\eta}}} = \int_{2}^{\infty}\sum_{1<q\le t}d(q)\,{{(\log t)+(1+\eta)}\over {t^2 (\log t)^{2+\eta}}}dt
 \ll_{\eta} \int_{2}^{\infty}{1\over {t (\log t)^{\eta}}}dt 
  \ll_{\eta} 1, 
$$
\par
\noindent
whence : bound $(2.1)$ above follows again from $\eta>1$.\hfill $\square$

\bigskip

\par
Theorem 2 entails the following result. 
\smallskip
\par
\noindent{\bf Proposition 2.} {\it The null-function, i.e., }\enspace  $\0:n\in \N \mapsto 0\in \C$ , {\stampatello may} $\underline{\hbox{\stampatello not}}$ {\stampatello have} {\it any} $\GRE$ {\stampatello with} {\it non-trivial coefficient, i.e. } \enspace $G\neq \0$, {\stampatello and $\etaDecay$ with } $\eta>1$. 
\smallskip
\par
\noindent{\bf Proof}. From Theorem 2, a $\GRE$ with $\eta>1$ for $\0$ {\stampatello entails} \enspace $G=\WinT \0=\CarT \0=\0$.\hfill $\square$

\vfill

\par
By the way, we have two non-trivial {\stampatello Ramanujan Coefficients for}\enspace $\0$, say
$$
R_0(q)\defineq {1\over q},
\enspace
\forall q\in \N
\quad
\hbox{\stampatello given by Ramanujan [R]}
\quad
\& 
\quad
H_0(q)\defineq {1\over {\varphi(q)}},
\enspace
\forall q\in \N
\quad
\hbox{\stampatello given by Hardy [H].}
$$

\bigskip

\par
A Good Ramanujan Expansion for $F$ entails also a non-trivial bound for the function expanded itself. 
\smallskip
\par
\noindent{\bf Theorem 3.} {\it Let } $F:\N \rightarrow \C$ {\it have a } {\stampatello Good Ramanujan Expansion}, {\it of } {\stampatello parameter } $\eta>0$. {\it Then}
$$
F(a)=F(1)+O_{\eta}\left(\sum_{{d|a}\atop {d>1}}{1\over {\log^{\eta} d}}\right),
\quad
\hbox{\stampatello uniformly}
\enspace \forall a\in \N. 
$$
\par
\noindent{\bf Remark 4.} We leave the Proof as an exercise for the diligent reader (hints in next Remark).\hfill $\diamond$
\par
\noindent{\bf Remark 5.} Actually, {\stampatello Th.m 3} follows from {\stampatello Th.m 1}, writing \enspace $\forall a\in \N$,\enspace  $F(a)=\sum_{d|a}F'(d)$.\hfill $\diamond$

\eject

\par				
The set of {\stampatello arithmetic functions} $F$ {\stampatello that divert values in } $\Z_{\sQ}^*$ , {\it for some fixed} $\sQ>2$, {\stampatello gives a} big {\stampatello set} of functions $F:\N \rightarrow \C$ {\stampatello without} $\GRE$, thanks to the following. 
\smallskip
\par
\noindent{\bf Theorem 4.} {\it Let } $F:\N \rightarrow \C$ {\stampatello divert values in } $\Z_{\sQ}^*$ , {\it for some} {\stampatello fixed} $\sQ>2$. {\it Then}
\smallskip
\item{$i)$} $\exists a\in \Z_{\sQ}^*\backslash \{1\}$ : $F(a)=F(p_a)\neq F(p_1)=F(1)$, \enspace $\forall p_a\equiv a(\!\bmod \,\,\sQ)$, $\forall p_1\equiv 1(\!\bmod \,\,\sQ)$ \enspace {\stampatello primes}
\smallskip
\item{$ii)$} {\stampatello the two subsequences}, {\it in the sequence of primes, } \thinspace $p_a\equiv a(\!\bmod \,\,\sQ)$ \thinspace {\stampatello and} \thinspace $p_1\equiv 1(\!\bmod \,\,\sQ)$ \enspace {\stampatello have}
$$
\lim_{p_a\equiv a(\!\!\bmod \sQ)}F'(p_a)\neq \lim_{p_1\equiv 1(\!\!\bmod \sQ)}F'(p_1)=0,
\enspace
\hbox{\stampatello whence}
\enspace
\lim_p F'(p),
\enspace 
\hbox{\stampatello over all primes, doesn't exist}
$$
\smallskip
\item{$iii)$} {\it in particular,} $F'$ {\stampatello is not infinitesimal,} {\it whence by Theorem 1}
\smallskip
\item{$iv)$} $F$ {\stampatello can't have any} $\GRE$.  

\medskip

\par
\noindent{\bf Proof}. Our \thinspace $i)$ \thinspace comes from the \lq \lq {\stampatello diverts values}\rq \rq-{\stampatello definition}.\hfill QED
\par
\noindent
The fact that there are infinitely many primes $p_a$ and $p_1$ in the above non-trivial arithmetic progressions modulo $\sQ$ comes straight from {\stampatello Dirichlet Theorem on primes in arithmetic progressions}, so, using above \thinspace $i)$ \thinspace and \enspace $F'(p)\EqByDef F(p)-F(1)$, $\forall p\in \Primes$, \enspace we get \thinspace $ii)$ \thinspace at once.\hfill QED
\par
\noindent
Finally, \thinspace $iii)$ \thinspace and \thinspace $iv)$ \thinspace {\stampatello are immediate}.\hfill $\square$
\medskip
\par
Like immediate are the applications of Theorem 4 to Correlations, following. 
\medskip

\bigskip

\par
\noindent{\bf 3. Application to Correlations} 
\bigskip
\par
\noindent
We wish to give our most interesting applications, namely for Correlations.
\medskip
\par
\noindent{\bf 3.1. A Correlation without the $\REEF$}
\medskip
\par
For general definitions about them, see [C1]. There (in second version), we define the $3-${\stampatello hypotheses correlations}, see $(0)$, $(1)$, $(2)$ hypotheses there; and, just to give, say, the easiest example of this kind, we build now what we call \lq \lq {\stampatello Counterexample} $1$\rq \rq.  This, in fact, is the Counterexample we need to prove that not all the $3-$hp.s Correlations have the {\stampatello R.E.E.F.}, abbrev.for {\stampatello Ramanujan Exact Explicit Formula}; this is : point $(iii)$ of Theorem 1 [CM] (this Theorem doesn't prove the REEF but gives four equivalent properties for the Ramanujan Expansion of a Correlation, $C_{f,g}(N,a)\defineq \sum_{n\le N}f(n)g(n+a)$, w.r.t. the {\it shift} $a\in \N$). Actually, the $\REEF$ is a fixed-length finite Ramanujan Expansion, entailing the Correlation having it to be a Truncated Divisor Sum (see [CM] considerations in sections 2,3 and 6, and compare [C0] Theorem 3, proving that fixed-length RE and TDS are the same,too). 
\par
Last but not least, following Counterexample 1 has been my source of inspiration for the \lq \lq {\stampatello diverting-values}\rq \rq \thinspace definition. In fact, it diverts values, whence has no G.R.E.; in particular, from this it's easy to prove (we leave it as an exercise for the careful reader) that no $\REEF$, for this Correlation, is possible. However, it's so easy that this last property will be checked at once. 
\medskip
\par
We start fixing a prime $p_0>2$, in correspondence of which we choose a length \thinspace $N\in \N$ \thinspace for our Correlation, with only one requirement : to be large enough, in order to allow the existence of another prime, {\stampatello say} \thinspace $n_0\in \Primes$, \thinspace $n_0>p_0$ \thinspace and \thinspace $n_0\equiv -1(\bmod \; p_0)$, \thinspace in the range \thinspace $n_0\le N$.
\medskip
\par
With these two primes and such a $N\in \N$, we are ready to build our {\stampatello Correlation}, say
$$
C_{f_0,g_0}(N,a)\defineq \sum_{n\le N}f_0(n)g_0(n+a), 
\qquad
\forall a\in \N,
$$
\par
\noindent
for arithmetic functions \thinspace $f_0\defineq \1_{\{n_0\}}$ \thinspace and \thinspace $g_0(m)\defineq c_{p_0}(m)$, \thinspace $\forall m\in \N$; getting
$$
C_{f_0,g_0}(N,a)\EqByDef c_{p_0}(a-1), 
\qquad
\forall a\in \N,
\leqno{(\hbox{\stampatello Counterexample } 1)}
$$
\par
\noindent
from \thinspace $n_0\equiv -1(\bmod \; p_0)$ \thinspace and \thinspace $p_0-${\stampatello periodicity of Ramanujan Sum} above. 
\par				
{\stampatello We may write}, by definition, \thinspace $\widehat{g_0}(q)=\1_{\{p_0\}}(q)$ \thinspace and \thinspace $g_0(m)=\sum_{q\le p_0}\widehat{g_0}(q)c_q(m)$, {\stampatello so we don't need to truncate divisors at } $Q\le N$, compare [C1], since our \thinspace $Q=p_0\le N$ \thinspace ($N$ is \lq \lq {\stampatello large enough}\rq \rq), {\stampatello and our Correlation is Fair}, trivially, {\stampatello because $a-$dependence is only inside Ramanujan Sum} above. Thus, hypothesis $(0)$, i.e., {\stampatello Basic Hypothesis} (see next $\S3.2$) is satisfied; the {\stampatello other two hypotheses} $(1)$ \& $(2)$ in [C1], namely $\widehat{g_0}$ is square-free supported for $(1)$ and hypothesis $(2)$, requiring \thinspace $f_0$ \thinspace supported on primes, are {\stampatello also both fulfilled}, from \thinspace $p_0,n_0\in \Primes$.\hfill In passing, \thinspace $g_0$ \thinspace {\stampatello is $p_0-$periodic, entailing} \thinspace $C_{f_0,g_0}$ \thinspace {\stampatello is $p_0-$periodic}, too. Compare next $\S3.2$ : in particular, {\it Wintner's Period} definition. 
\smallskip
\par
This is a kind of {\stampatello toy-model Correlation}, so its {\stampatello periodicity is very easy}. 
\par
\noindent
From [C1], see section $4$ after $(\ast)$ there, we know that for general, \lq \lq {\stampatello real-world sounding}\rq \rq \thinspace $3-$hypotheses Correlations, we still have periodicity, linked to the {\stampatello truncation parameter} $Q\in \N$, with {\stampatello period} {\tt T}$\in \N$ ({\stampatello even more difficult to find}). Compare quoted $\S3.2$. 
\smallskip
\par
We'll come to a general approach to {\stampatello periodicity of very general correlations}, namely only having Basic Hypothesis, in future publications regarding Correlations. Again, compare $\S3.2$.

\bigskip

\par
\noindent
Well, \thinspace $C_{f_0,g_0}$ \thinspace {\stampatello diverts values in} \thinspace $\Z_{p_0}^*$. And, being very simple, it does it {\stampatello in the simplest way}, namely: 
\smallskip
\par
for \qquad $a=1\in \Z_{p_0}^*$,\enspace $C_{f_0,g_0}(N,a)=c_{p_0}(a-1)=\varphi(p_0)=p_0-1$ 
\smallskip
and \qquad $\forall a\in \Z_{p_0}^*\backslash \{1\}$,\enspace $C_{f_0,g_0}(N,a)=c_{p_0}(a-1)=\mu(p_0)=-1$. 
\medskip
\par
\noindent
Thus Theorem 4 tells us that this Correlation has {\stampatello Eratosthenes Transform}
$$
C'_{f_0,g_0}(N,d)\defineq \sum_{a|d}C_{f_0,g_0}(N,a)\mu\left({d\over a}\right)
 =\sum_{a|d}c_{p_0}(a-1)\mu\left({d\over a}\right)
\qquad
\hbox{{\stampatello not infinitesimal}}
\enspace 
(\hbox{\it as}\enspace d\to \infty)
$$
\par
\noindent
and, in particular, it has no G.R.E. (from Theorem 1).\hfill Next considerations exclude the $\REEF$, for it. 

\bigskip

\par
In fact, any Basic Hypothesis Correlation, say, \thinspace $C_{f,g_Q}(N,a)$,\thinspace with the 
$$
C_{f,g_Q}(N,a)=\sum_{q\le Q}{{\widehat{g_Q}(q)}\over {\varphi(q)}}\sum_{n\le N}f(n)c_q(n)c_q(a),
\qquad
\forall a\in \N,
\leqno{\REEF}
$$
\par
\noindent
compare [C1] Theorem 3.1 for the coefficients, has clearly {\stampatello infinitesimal Eratosthenes Transform}, since {\stampatello the $\REEF$ entails} (another exercise for the careful reader, applying $(1.7)$ in [C1]: Kluyver Formula)
$$
C'_{f,g_Q}(N,d)=0,
\qquad
\forall d>Q. 
$$
\par
\noindent
The same result quoted above allows to calculate our Correlation's coefficients; say, Carmichael-Wintner $\ell-$th coeff.'s in general:
$$
{{\widehat{g_Q}(\ell)}\over {\varphi(\ell)}}\sum_{n\le N}f(n)c_{\ell}(n)
$$
\par
\noindent
as above; {\stampatello for our case $f=f_0$ \& $g_Q=g_0$, it's}
$$
\1_{\{p_0\}}(\ell)\cdot {1\over {\varphi(\ell)}}c_{\ell}(n_0)
$$
\par
\noindent
giving a $\REEF$'s RHS($=$Right Hand Side), for our Correlation, 
$$
\sum_{\ell \le p_0}\1_{\{p_0\}}(\ell)\cdot {1\over {\varphi(\ell)}}c_{\ell}(n_0)c_{\ell}(a)
=
{{\mu(p_0)}\over {\varphi(p_0)}}c_{p_0}(a).
$$
\par
\noindent
Thus
$$
C_{f_0,g_0}(N,a)=c_{p_0}(a-1)
 =\sum_{\ell \le p_0}\1_{\{p_0\}}(\ell)\cdot {1\over {\varphi(\ell)}}c_{\ell}(n_0)c_{\ell}(a)
  ={{\mu(p_0)}\over {\varphi(p_0)}}c_{p_0}(a)
\qquad
\hbox{\stampatello can't hold}
\enspace
\forall a\in \N. 
$$

\medskip

\par				
\noindent
Namely,

\medskip

\par
{\stampatello our Correlation \thinspace $c_{p_0}(a-1)$ \thinspace has no $\REEF$, in perfect agreement with Theorem 4}.

\medskip

\par
\noindent
In fact, RHS above is
$$
-{1\over {p_0-1}}c_{p_0}(a)={1\over {p_0-1}},
\qquad
\forall a\not\equiv 0(\bmod \; p_0), 
$$
\par
\noindent
while our Correlation's absolute value is
$$
\left|c_{p_0}(a-1)\right|\ge 1,
\qquad
\forall a\in \N, 
$$
\par
\noindent
since {\stampatello modulus $p_0$ is square-free} (use $(1.6)$ in [C1]: H\"older's Formula). 

\medskip

\par
Just to have fun, say, the \lq \lq $\REEF$ {\stampatello holds} $\forall a\in \N,a\equiv 0(\bmod \; p_0)$\rq \rq. A kind of \lq \lq {\stampatello Relative REEF}\rq \rq. 

\bigskip

\par
\noindent
More in general, {\stampatello our Theorem 4 entails not only that diverting values $C_{f,g_Q}(N,a)$ have no $\REEF$; but also that their IPPification $\widetilde{C_{f,g_Q}}(N,a)$} (compare [C1]), {\stampatello with same Eratosthenes Transform, but supported on square-free divisors} (see [C1],[C2]), {\stampatello is still diverting-values, so no $\REEF$, for $\widetilde{C_{f,g_Q}}(N,a)$ too}. 

\medskip

\par
Again, this is immediate for our Toy-Model Correlation (as the careful reader may check, recalling IPPification definition in [C1],[C2]). 
\medskip

\bigskip

\par
\noindent{\bf 3.2. Any $\BH$ Correlation with a $\GRE$ has the $\REEF$} 
\bigskip
\par
\noindent
We {\stampatello recall} that a $\BH$ Correlation $C_{f,g_Q}(N,a)$ satisfies the hypothesis $(0)$ quoted above, i.e.: 
\smallskip
\par
\centerline{BASIC HYPOTHESIS $\BH$}
\smallskip
\par
\noindent
$g_Q$ {\stampatello has} Eratosthenes Transform $g_Q'$ with $\supporto(g_Q')\subseteq [1,Q]$ : $Q\in \N$ is the {\stampatello range} of $g_Q$, {\stampatello satisfying} $Q\le N$
\smallskip
\par
\centerline{\bf and}
\smallskip
\par
\noindent
$$
 C_{f,g_Q}(N,a)=\sum_{q\le Q}\widehat{g_Q}(q)\sum_{n\le N}f(n)c_q(n+a),
\quad 
\forall a\in \N,
$$
\par
\noindent
in which RHS {\it depends on $a\in \N$ only inside } $c_q(n+a)$ (equivalently, nor $\widehat{g_Q},f$, neither their supports depend on $a\in \N$ in above equation) : {\stampatello Correlation is Fair} 
\par
\noindent
(Compare [C1] to deepen why {\stampatello range up to $Q\le N$ and Fairness} are V.I.P., Very Important Properties) 
\smallskip
\par
A $\BH$ Correlation has the Ramanujan Exact Explicit Formula IFF (by definition) page 7's $\REEF$ holds. 
\smallskip
\par
We {\bf define} (compare [C1] definitions) {\stampatello Wintner's Period} for all $g_Q\neq \0$ (avoiding $\widehat{\0}=\0$) as above : 
$$
W(g_Q)\defineq {\rm l.c.m.}\{q\le Q : \widehat{g_Q}(q)\neq 0\}
$$ 
\par
\noindent
while in [C1] (see version 2) we introduced {\script Q}$\defineq {\rm l.c.m.}\{q\le Q\}$ in Theorem 3.1 Proof : see $\S4.1$ there. 
\medskip
\par
By definition, \enspace $g_Q(m)=\sum_{q\le Q}\widehat{g_Q}(q)c_q(m)$, $\forall m\in \N$\enspace [C1]\enspace is $W(g_Q)-${\stampatello periodic}; however, {\stampatello equation above entails} $C_{f,g_Q}(N,a+W(g_Q))=C_{f,g_Q}(N,a)$, $\forall a\in \N$: our {\stampatello Correlation's $W(g_Q)-$periodic, too.}

\vfill

\par
As a Corollary to Corollary 2 above, for $F(a):=C_{f,g_Q}(N,a)$, $\forall a\in \N$, we get the following {\stampatello Application to Correlations} giving the title to present subsection. We use: this $F$ is $r-${\stampatello periodic}, with $r=W(g_Q)$; so, $q|r \Rightarrow q\le Q$, in $(\ast)$ above. In fact, $\widehat{g_Q}(q)=0$ entails $\Win_q \,F=0$ (again, same $F$), from Theor.3.1 [C1]. 
\smallskip
\par
\noindent{\bf Proposition 3.} {\it Any $\BH$ Correlation with a $\GRE$ has the $\REEF$.} 
\medskip
\par
Actually (but we leave the Proof as an exercise to the careful Reader), even more: 
\smallskip
\par
\noindent{\bf Proposition 4.} {\it Let $C_{f,g_Q}(N,a)$ satisfy $\BH$. Then,}
$$
\hbox{\stampatello the eventual } \GRE \enspace \underline{\hbox{\stampatello is}} \enspace \hbox{\stampatello the eventual \REEF.} 
$$

\eject

\par				
\noindent
Previous subsection shows, {\stampatello through Counterexample $1$}, that {\stampatello a $\BH-$Correlation might be so irregular to divert values}, whence no $\GRE$. On the other hand, we see from Propositions 3 \& 4 that {\stampatello if our $\BH-$Correlation}, so to speak, {\stampatello has any kind of} $\GRE$, i.e., with any $\eta>0$ fixed (otherwise $\eta>1$ from Theorem 2 entails the $\REEF$, rather at once), this becomes, \lq \lq {\stampatello at once}\rq \rq, {\stampatello the $\REEF$}; testifying, in this way, that our Correlation has a {\stampatello finite Ramanujan expansion} (actually, of {\stampatello fixed-length}) with Carmichael-Wintner Coefficients. 
\medskip
\par
In other words, there's {\stampatello a kind of} \lq \lq polarization\rq \rq, with two, very different cases for a $\BH-$Correlation.  
\medskip
\par
\noindent
We {\stampatello may ask} for a \lq \lq minimum level\rq \rq, of {\stampatello regularity, i.e. a $\GRE$}; but this is {\stampatello possible only with the $\REEF$}; and, on the converse, {\stampatello as soon as there's no $\REEF$, no $\GRE$ can be found.}
\medskip
\par
In particular, since {\stampatello having $q-$th Ramanujan Coefficients in a shift-Ramanujan Expansion}: (compare [CM] terminology and results)
$$
C_{f,g_Q}(N,a)=\sum_{q=1}^{\infty}R(q)c_q(a),
\qquad
\forall a\in \N
\leqno{\hbox{\stampatello (SRE)}_{\delta}}
$$
\par
\noindent
{\stampatello in the so-called \lq \lq $\delta-$class\rq \rq}, for a {\stampatello fixed } $\delta>0$, namely {\stampatello requiring}
$$
R(q)\ll_{\delta} {1\over {q^{1+\delta}}}, 
\qquad
\hbox{\stampatello as}\enspace q\to \infty
\leqno{(\hbox{$\delta-${\stampatello class}})}
$$
\par
\noindent
{\stampatello entails} a $\GRE$ with $G=R$ \enspace as a {\stampatello Ramanujan Coefficient, implying} together with Proposition 3: 
$$
\hbox{\stampatello (SRE)}_{\delta}
\enspace \Longrightarrow \enspace
\GRE
\enspace \Longrightarrow \enspace
\REEF.
$$
\par
So, in [CM] {\stampatello we} still {\stampatello didn't know that the existence of Ramanujan Coefficients in $\delta-$class} (see above bound) {\stampatello for the Shift-Ramanujan Expansion} (above {\stampatello (SRE)}$_{\delta}$ equation), actually, {\stampatello entails the $\REEF$.} (Actually, this name came out later: see $(iii)$ in our [CM] Theorem 1. We recall it abbreviates {\stampatello Ramanujan Exact Explicit Formula}.) 
\medskip

\bigskip

\par
\noindent{\bf 4. Final Remarks, further properties and coming soon} 
\par
\noindent
\bigskip
We start with a Remark about how this paper was born. 
\smallskip
\par
\noindent{\bf Remark 6.} Regarding the Correlations, in our paper with Ram Murty [CM] about the finite Ramanujan expansions, we assumed that the coefficients in their {\it shift-Ramanujan expansion} had an extremely fast decay, gaining, as we assumed, powers of the index $q$ (NOT $\log-$POWERS, like here). Then, playing with this assumption, I was wondering the effect, over an arbitrary arithmetic function, of having a much slower decay, namely the $\etaDecay$ above. Subsequently, I realized Theorem 1, namely that it descends, from this decay for the Ramanujan Coefficients $G(q)$ of our $F$, that $F'(d)$ is infinitesimal as $d\to \infty$, compare Th.m 1. In particular, $F'(p)\EqByDef F(p)-F(1)$ is infinitesimal as $p\to \infty$, along primes $p$; but, playing around with formul\ae, for ($3-$Hypotheses-)Correlations in [C1], I saw a high probability to get oscillating such $F'(p)$ along primes $p\to \infty$; last but not least, again from [C1], I knew that a Basic Hypothesis Correlation is $\sQ-$periodic, see [C1] : looking at residue classes in $\Z_{\sQ}^*$ and applying Theorem 4 arguments, I got the non-existence, in diverting-values case, of $F'(p)$ limit. In fact, $F'(p)$ was {\bf diverting values} of limits for it, in the subsequences for different residue classes ($1$ and $a\neq 1$, if $F(a)\neq F(1)$) : Thanks very much, to our Peter Gustav Lejeune Dirichlet ! Otherwise, no subsequences !\hfill Actually, these are vital for our present Corollary 2 Proof. 
\par
\noindent
In fact, along them we calculate the limit of our infinitesimal Eratosthenes Transform, getting from this hypothesis that \lq \lq $F$ is $r-$periodic\rq \rq, so to speak, becomes \lq \lq $F$ is $r-$even\rq \rq.\hfill See above Corollary 2 Proof.\hfill $\diamond$ 
\medskip
\par
By the way, as it's clear from Theorem 4 above, {\bf in case the Correlation diverts values} in $\Z_{\sQ}^*$, we may forget about the Possibility, not only to have it in First Class (see [CM]), but also, of course, to get a Good Ramanujan Expansion, for it ! However, see next Remarks 10 \& 11. 
\medskip
\par
We get back, by the way, to the $\etaDecay$, since we get an interesting consequence of Theorem 2. 
\smallskip
\par				
\noindent{\bf Remark 7.} From Th.m 2, we have decay for the $q-$th Wintner Coefficient, now $\Win_q\, F=G(q)$, that is, from $\etaDecay$,\enspace $\Win_q\, F\ll_{\eta} (\log q)^{-1-\eta}/q$, of course; however, since Th.m 1 gives: \enspace  $F'(d)\ll_{\eta} (\log d)^{-\eta}$, the series defining \enspace $\Win_q\, F$ \enspace has the upper bound, for all $q\ge 2$, (compare above Theorems' Proofs' calculations)
$$
|\Win_q\, F|\EqByDef {1\over q}\left|\sum_{m=1}^{\infty}{{F'(qm)}\over m}\right|
 \le {1\over q}\sum_{m=1}^{\infty}{{\left|F'(qm)\right|}\over m}
  \ll_{\eta} {1\over q}\,\sum_{m\le q}{1\over {m(\log q)^{\eta}}} + {1\over q}\,\sum_{m>q}{1\over {m(\log m)^{\eta}}}
   \ll_{\eta} (\log q)^{1-\eta}/q, 
$$
\par
\noindent
easily getting a {\stampatello $\log-$squared gain for the decay}, motivated by the sign-change (in case $F:\N \rightarrow \R$) or, more generally, phase-change of our Eratosthenes Transform, within these series. Say, our Th.m 2 hypotheses \lq \lq {\stampatello force the cancellation}\rq \rq, coming from \thinspace $F'(d)$, in these series defining Wintner coefficients.\hfill $\diamond$ 
\medskip
\par
We continue with two remarks about the $\GRE$ definition. 
\smallskip
\par
\noindent{\bf Remark 8.} From the bound \enspace $|c_q(a)|\le (a,q)$ (see Lemma A.1 in [CM]), we get, applying $\etaDecay$, 
$$
\forall a\in \N,
\quad
\sum_{q=1}^{\infty}\left|G(q)c_q(a)\right|\le \sum_{q=1}^{\infty}(a,q)\left|G(q)\right|
 \ll_{\eta} 1+d(a), 
$$
\par
\noindent
whence we get {\stampatello the Absolute Convergence} of our $\GRE$. Then, {\stampatello in particular, all summation methods for the $\RamaExp$ give} $\underline{\hbox{\stampatello both}}$ {\stampatello convergence and to the same sum}.\hfill $\diamond$ 
\medskip
\par
\noindent{\bf Remark 9.} In the $\GRE$ definition we have included $\LuchtExp$, because it's not immediately descending from $\RamaExp$, as {\stampatello applying Kluyver Formula} (quoted $(1.7)$ in [C1])
$$
F(a)=\sum_{q=1}^{\infty}G(q)c_q(a)
 =\sum_{d|a}d\sum_{{q=1}\atop {q\equiv 0(\!\!\bmod d)}}^{\infty}G(q)\mu\left({q\over d}\right)
  =\sum_{d|a}d\sum_{K=1}^{\infty}\mu(K)G(dK) 
$$ 
\par
\noindent
{\stampatello doesn't imply} $\LuchtExp$ : {\stampatello in case} $G$ {\stampatello depends on $a$, we can't apply M\"obius Inversion}.\hfill $\diamond$
\medskip
\par
We come now to other two remarks, pertaining to Correlations. 
\smallskip
\par
\noindent{\bf Remark 10.} For {\stampatello Counterexample} 1, this Correlation is a {\stampatello shifted Ramanujan Sum}, see $\S3$, giving at once the {\stampatello diverting-values property}; for, say, \lq \lq real-life Correlations\rq \rq, the check of this property is {\stampatello not easy \& may require some new hypotheses}.\hfill $\diamond$
\medskip
\par
\noindent{\bf Remark 11.} The {\stampatello Hardy-Littlewood Correlation} (for $a=2k-$twin primes, say): 
$$
C_{\Lambda,\Lambda}(N,a)\defineq \sum_{n\le N}\Lambda(n)\Lambda(n+a),
\quad
\forall a\in \N, 
\qquad
\hbox{\stampatello has} 
\qquad
C_{\Lambda,\Lambda_N}(N,a)\defineq \sum_{n\le N}\Lambda(n)\Lambda_N(n+a),
\quad
\forall a\in \N, 
$$
\par
\noindent
as its \lq \lq {\stampatello truncation}\rq \rq \thinspace [C1], where $\Lambda$ is the von-Mangoldt function (see [C1]); this truncation is {\stampatello easier to study}, compare Theorem 3.1 [C1], but the check to verify the {\stampatello Monochromaticity} is not easy at all. However, this problem is very interesting, thanks to Theorem 4, as in case of a {\stampatello diverting-values Correlation}, we know that (see $\S3$ considerations, above) not only, say, $C_{f,g_Q}(N,a)$ {\stampatello may not have the $\REEF$}, but also $\widetilde{C_{f,g_Q}}(N,a)$ {\stampatello has no} $\REEF$.\hfill $\diamond$
\medskip
\par
\noindent
We end up with Arithmetic Functions: any $\RamaExp$ has no coefficients-decay ensured. 
\smallskip
\par
\noindent{\bf Remark 12.} Once {\stampatello fixed} $F:\N \rightarrow \C$, a result of Hildebrand dated 1984 finds {\stampatello a} $\RamaExp$ for $F$; actually, a {\stampatello finite} one [C2],[CM], {\stampatello but not a fixed-length one} [C0] : in particular, no Carmichael-Wintner Coefficients; a complete exposition, of this result \& full details, is in the Appendix of [C3].\hfill $\diamond$
\medskip
Before a short deepening of Proposition 1 results, we point out: nowadays definition, in $\S1$, of \lq \lq $r-${\stampatello even}\rq \rq, actually, is the one given by Eckford Cohen [Coh] as $(3)$, compare His paper; also, the Ramanujan Expansion we apply \& prove with \lq \lq our methods\rq \rq, actually is recognized as equivalent to quoted $(3)$ and labeled as $(1)$ in quoted paper. Actually, we recognized this Expansion's Coefficients as Carmichael-Wintner ones. In [Coh], for example $(2)$ \& other properties are proved to be equivalent to the nowadays definition of $r-${\stampatello even}.

\par				
\noindent
We wish, here, to deepen Corollary 2 and, also, Proposition 1 : with a more general requirement for $F$ than assuming \enspace $F'(D)\to 0$, {\stampatello as $D\to \infty$ in natural numbers}. We consider the constraint on $F$ for a particular class of natural numbers $D$, going to infinity. We'll see soon: once assumed $F$ is $r-${\stampatello periodic}, the following 
$$
\forall d|r, \enspace \lim_p F'(dp)=0
\leqno\Cdp
$$ 
\par
\noindent
({\stampatello with } $\lim$ \enspace {\stampatello over }$p\in \Primes$, $p\to \infty$) is equivalent to : $F$ is $r-${\stampatello even}. We avoid trivial cases $r=1$ \& $r=2$. 
\smallskip
\par
\noindent{\bf Theorem 5.} {\it Let } $r\in \N$, $r>2$ {\it and let } $F:\N \rightarrow \C$ {\it be } $r-${\stampatello periodic}. Then, 
\smallskip
\par
\centerline{$F$ {\it satisfies} $\Cdp$ \enspace $\Longleftrightarrow$ \enspace $F$ {\it is } $r-${\stampatello even}.} 
\smallskip
\par
\noindent{\bf Proof}. Our implication \lq \lq $\Rightarrow$\rq \rq \thinspace has the same Proof working for Corollary 2; however, now we give details. In $\Cdp$ the prime $p\to \infty$ in all $\Primes$, while in quoted Proof it's in the subsequence $p\equiv c\,(\bmod \; r/d_0)$: in all subsequences, the limit exists and is the same. We profit to clarify that, compare Corollary 1 Proof,
$$
p>d_0
\enspace \Rightarrow \enspace 
F(d_0p)\EqByDef\sum_{t|(d_0p)}F'(t)
        =\sum_{d|d_0}\left(F'(dp)+F'(d)\right)
         =\sum_{d|d_0}F'(dp)+F(d_0). 
$$
\vskip-0.54cm
\par
\hfill QED
\par
\noindent
For \lq \lq $\Leftarrow$\rq \rq: $\forall a\in \N$, $F(a)=F((a,r))\EqByDef\sum_{t|a,t|r}F'(t)$, whence we have $F'(t)=0$, $\forall t>r$; in turn, this  entails $F'(dp)=0$, $\forall d\in \N$, $\forall p>r$.\hfill $\square$ 
\medskip
\par
Notice that $\Cdp$, for the case $d=1$, requires \enspace $\lim_p F'(p)=0$, {\stampatello i.e.}, from $F'(p)\EqByDef F(p)-F(1)$, $\forall p\in \Primes$, $\forall F:\N \rightarrow \C$, \enspace $\lim_p F(p)=F(1)$; whence, with $r-${\stampatello periodicity} entails: $F$ is $r-${\stampatello Monochromatic}. 
\par
This is not a great surprise, from Theorem 4 results. However, above Condition has still all the other $d|r$ \& $d>1$ requirements on our $F$ !!! 
\par
\noindent
(So, without additional Hypotheses we expect that $F$ being $r-${\stampatello Monochromatic} is a weaker property than $F$ is $r-${\stampatello even}; recall that saying $F$ is $r-${\stampatello Monochromatic} assumes $F$ is $r-${\stampatello periodic}.)  
\medskip
\par
From Remarks above and the results in this paper, especially present Corollary 2, it is evident that our approach is still very fertile, even in very general results: compare Proposition 1 \& Theorem 5, for very general arithmetic functions; however, {\stampatello diverting-values Correlations} $\underline{\hbox{\stampatello and}}$ {\stampatello the $\sQ-$Monochromatic ones} are still worthy to study ! 
\medskip
\par
This we will pursue in future papers. 
\medskip
\vfill
\par
\noindent{\bf Acknowledgment.} I profit, here, to thank once again Ram Murty, starting from his survey [M] introducing me to the Ramanujan Expansions; and for the patience, in working with me [CM] on Finite Ramanujan Expansions applied to Correlations, namely shifted convolution sums. And this [CM], also, was an introduction for me, to the world of Ramanujan Expansions \& Correlations ! 

\eject

\par				
\centerline{\stampatello Bibliography}
\medskip
\item{[Coh]} E. Cohen, {\sl A class of arithmetical functions}, Proc. Nat. Acad. Sci. U.S.A. {\bf 41} (1955), 939--944.
\smallskip
\item{[C0]} G. Coppola, {\sl Recent results on Ramanujan expansions with applications to correlations}, Rend. Sem. Mat. Univ. Pol. Torino {\bf 78.1} (2020), 57--82. 
\smallskip
\item{[C1]} G. Coppola, {\sl General elementary methods meeting elementary properties \thinspace of \thinspace correlations}, {\tt available}\break{\tt online at} \enspace arXiv:2309.17101 (2nd version)
\smallskip
\item{[C2]} G. Coppola, {\sl On Ramanujan smooth expansions for a general arithmetic function}, {\tt available online at} \enspace arXiv:2407.19759v1 
\smallskip
\item{[C3]} G. Coppola, {\sl Absolute convergence of Ramanujan expansions admits coefficients' coexistence of Ramanujan expansions}, {\tt available online at} \enspace arXiv:2502.14415v1 
\smallskip
\item{[CM]} G. Coppola and M. Ram Murty, {\sl Finite Ramanujan expansions and shifted convolution sums of arithmetical functions, II}, J. Number Theory {\bf 185} (2018), 16--47. 
\smallskip
\item{[H]} G.H. Hardy, {\sl Note on Ramanujan's trigonometrical function $c_q(n)$ and certain series of arithmetical functions}, Proc. Cambridge Phil. Soc. {\bf 20} (1921), 263--271.
\smallskip
\item{[HL]} G.H. Hardy and J.E. Littlewood, {\sl SOME PROBLEMS OF 'PARTITIO NUMERORUM'; III: ON THE EXPRESSION OF A NUMBER AS A SUM OF PRIMES.} Acta Mathematica {\bf 44} (1923), 1--70. 
\smallskip
\item{[M]} M. R. Murty, {\sl Ramanujan series for arithmetical functions}, Hardy-Ramanujan J. {\bf 36} (2013), 21--33. Available online 
\smallskip
\item{[R]} S. Ramanujan, {\sl On certain trigonometrical sums and their application to the theory of numbers}, Transactions Cambr. Phil. Soc. {\bf 22} (1918), 259--276.
\smallskip
\item{[T]} G. Tenenbaum, {\sl Introduction to Analytic and Probabilistic Number Theory}, Cambridge Studies in Advanced Mathematics, {46}, Cambridge University Press, 1995. 
\medskip
\vfill
\par
\leftline{\tt Giovanni Coppola - Universit\`{a} degli Studi di Salerno (affiliation)}
\leftline{\tt Home address : Via Partenio 12 - 83100, Avellino (AV) - ITALY}
\leftline{\tt e-mail : giocop70@gmail.com}
\leftline{\tt e-page : www.giovannicoppola.name}
\leftline{\tt e-site : www.researchgate.net}
\eject
\bye